\documentclass[a4paper,12pt]{article}
\usepackage{amsmath}
\usepackage{amsthm}
\usepackage{epsfig}
\usepackage{psfrag}
\usepackage{subfigure}
\usepackage[mathscr,mathcal]{eucal}
\usepackage{amssymb}
\usepackage{enumerate}
\def\phi{\varphi}

\def\be{\begin{equation}}
\def\ee{\end{equation}}

\def\bea{\begin{eqnarray}}
\def\eea{\end{eqnarray}}

\def\cO{\mathcal{O}}

\def\cO{\mathcal{O}}

\def\diff{\,\mathrm{d}}

\def\sigmal{\underline{\sigma}}
\def\sigmau{\overline{\sigma}}

\def \R {\mathbb R}

\def \ind{1\!\!1}

\newcommand{\prob}[1]{\ensuremath{\mathbf{P}\big(#1\big)}}    
    
\newcommand{\var}[1]{\ensuremath{\mathbf{Var}\big(#1\big)}}

\def\usigma{\underline{\sigma}}
\def\osigma{\overline{\sigma}}

\def\vareps{\varepsilon}

\def\I{I}
\def\II{I\!I}

\def\Opm{\Omega^{\pm}}
\def\OI{\Omega^{\I}}
\def\OII{\Omega^{\II}}

\def\op{\omega^{+}}
\def\om{\omega^{-}}
\def\opm{\omega^{\pm}}
\def\oI{\omega^{\I}}
\def\oII{\omega^{\II}}

\def\mp{\mu^{+}}
\def\mm{\mu^{-}}
\def\mpm{\mu^{\pm}}
\def\mIM{\mu^{\I,M}}
\def\mIIc{\mu^{\II,c}}
\def\mII{\mu^{\II}}

\def\SIM{{\cal S}^{\I,M}}
\def\SIIc{{\cal S}^{\II,c}}

\def\VIM{V^{\I,M}}
\def\QIM{Q^{\I,M}}
\def\VIIc{V^{\II,c}}

\def\QIIc{Q^{\II,c}}

\def\tVIM{\widetilde V^{\I,M}}
\def\tQIM{\widetilde Q^{\I,M}}
\def\QII{Q^{\II}}

\newtheorem {theorem}{Theorem}

\newtheorem* {theorem*}{Theorem}
\newtheorem* {thm*}{Theorem}
\newtheorem* {lemma*}{Lemma}
\newtheorem* {corollary*}{Corollary}
\newtheorem* {prop*}{Proposition}
\newtheorem* {definition*}{Definition}
\newtheorem* {remark*}{Remark}

\begin{document}

\title{On the zero mass limit of tagged particle \\
diffusion in the 1-d Rayleigh-gas}
\author{
\sc P\'eter B\'alint$^\text{1,4}$ \qquad 
B\'alint T\'oth$^\text{1}$ \qquad 
P\'eter T\'oth$^\text{2,3}$
\\[5pt]
$^\text{1}$Institute of Mathematics (IM), Technical University of
  Budapest (BME)
\\[3pt]
$^\text{2}$R\'enyi Institute (RI) of the  Hungarian Academy of Sciences (MTA)
\\[3pt]
$^\text{3}$ Stochastics Research Group of MTA  at  IM BME
\\[3pt]
$^\text{4}$Courant Institute, New York University} 

\date{\today}

\maketitle

\centerline
{\large\sl Dedicated to Domokos Sz\'asz on his 65th birthday}

\bigskip

\begin{abstract}
We consider the $M\to0$ limit for tagged particle diffusion in a
1-dimensional Rayleigh-gas, studied originaly by Sinai and
Soloveichik \cite{sinaisoloveichik}, respectively  by Sz\'asz and
T\'oth  \cite{szasztoth1}. In this limit we derive
a new type of 
model for tagged paricle diffusion, with Calogero-Moser-Sutherland
(i.e. inverse quadratic) interaction 
potential between the two central particles. Computer simulations on
this new model reproduce exactly the numerical value 
of the limiting variance obtained by Boldrighini, Frigio and Tognetti
in \cite{boldrighinifrigiotognetti}.  
\end{abstract}

\section{Introduction}
\label{s:introduction}

The problem
of deriving relevant information on the diffusive scaling limit of
tagged particle motion (i.e. self-diffusion), from microscopic
principles has been
undoubtedly at the heart of mathematically rigorous statistical
physics of time dependent phenomena, at least since Einstein's
groundbreaking work.  
Mathematically rigorous investigation of tagged particle diffusion in
systems of particles governed by deterministic (Hamiltonian) dynamics
is notoriously difficult even in one dimensional models. After
remarkable advances achieved up to the late 1980-s (see
Section~\ref{s:oldresults} below and references cited there), in the
last twenty years there seemed to be less intense activity in the
field. This is certainly due to the difficulty of these problems and
lack of technical tools to attack them. 

In the present note we make a small  but hopefully not completely
irrelevant contribution to the subject. We investigate the $M\to0$
small mass limit of the tagged particle diffusion in the so-called
1-dimensional Rayleigh-gas. This system consists of an infinitely
extended one-dimensional system of point-like particles of mass 1 and
one single tagged particle of mass $M$ immersed in it. The particles
perform uniform motion and interact through elastic collisions. The
system is distributed according to the equilibrium Gibbs
measure. This means independent exponentially distributed
inter-particle distances and independent normally distributed
velocities with mean zero and inverse mass variances (that is:
equidistributed kinetic energy). It  is a fact that the dynamics of
the infinitely extended system is almost surely well defined under this
stationary measure. That is: no
multiple collisions and no accumulation of infinitely many particles
in finite time occurs. Randomness comes into the problem only through
the thermal equilibrium of the initial condition, otherwise the
dynamical evolution is deterministic. The central question is
understanding the diffusive scaling limit of the trajectory of the
tagged particle: $A^{-1/2}Q_{At}$, as $A\to\infty$. There exist a
number of deep and interesting results related to this problem which
will be shortly surveyed in Section~\ref{s:oldresults}. In the present
note we investigate the limit when $M\ll1$. We prove that in this
limit the system becomes equivalent to another, new model of tagged
particle motion, differing from the one described above in having
instead of one central particle of different mass, all particles of
the same mass but the
two central particles interacting via a Calogero-Moser-Sutherland-type
repulsive potential, with random strength parameter.  This result
explains some phenomena observed in earlier computer simulations on
the Rayleigh-gas. In particular the instability observed for small
values of $M$. We also present numerical simulations on this new
model. Our simulation results reproduce very accurately the numerical
value of the limiting variance of tagged particle in the Rayleigh-gas,
in the $M\to0$ limit, which were obtained in
\cite{boldrighinifrigiotognetti}. We claim that our result not only
reproduces (from a completely different approach) the numerical value
but also gives theoretical explanation of the phenomenon. 

The paper is organised as follows: 
In Section~\ref{s:models} we define the models of interacting particle
systems considered, their stationary Gibbs measures and the stochastic
processes whose diffusive asymptotics is later analysed. In
Section~\ref{s:oldresults} we briefly survey the existing earlier
results (rigorously proved and numerical, alike) on tagged particle
diffusion in 1-d Rayleigh-gas. In Section~\ref{s:Mtozero} we properly
state and prove the theorem which states that in the $M\to0$ limit the
1-d dynamics of the Rayleigh-gas with tagged particle of mass $M$
converges (trajectory-wise, in a natural topology) to the dynamics of
the 1-d gas of 
particles with Calogero-Moser-Sutherland interaction between the two
central particles. Finally, in Section~\ref{s:newnumerics} we present
our new numerical results referring to this new type of interacting
particle system. We should emphasize here that our numerical results
are not just accurate reproduction of older computer experiments but
are performed on a genuinely new type of model. One of the main points
of this paper is exactly the fact that these genuinely new numerical
results are in accurate agreement with the results of Boldrighini,
Frigio, Tognetti \cite{boldrighinifrigiotognetti}, giving independent
enhancement \emph{and theoretical explanation} to them.

\section{Models: state space, dynamics, stationary measures}
\label{s:models}

In this section we describe the models considered throughout the paper. In 
section~\ref{ss:statespace} we present a formal definition of the state 
spaces and the natural measures on them. Section~\ref{ss:dynamics}, which gives 
a more verbal description of the time evolution in our dynamical systems, clarifies 
that these models indeed correspond to the one dimensional gases mentioned 
in the Introduction.

\subsection{State spaces and stationary Gibbs measures}
\label{ss:statespace}

Let
\bea
\notag
\Opm :=
\{ \opm = (x_{\pm i},v_{\pm i})_{i=1}^{\infty}: 
(x_{\pm i},v_{\pm i})\in\R^\pm\times\R, \ \ 
x_{\pm1}=0,  \ \ 
\pm(x_{\pm(i+1)}-x_{\pm i})\ge0 \}. 
\eea
With slight abuse of notation and terminology sometimes we don't
distinguish between $\opm$ and the (unordered) set of points
$\{(x_{\pm i},v_{\pm i}): i=1,2,\dots\}\subset\R^\pm\times\R$. 
We endow the spaces $\Opm$ with the
topology defined by pointwise convergence: $\opm_n\to\opm$
iff $(x_{\pm i},v_{\pm i})_n\to(x_{\pm i},v_{\pm i})$, for all
$i=1,2,\dots$.  
This is a metrizable topology and
makes $\Opm$ complete and separable (i.e. Polish) spaces. 

We denote by $\mpm$ the following probability measures over
$\Opm$, respectively. Under $\mpm$ the random variables
$\xi_{\pm i}:=\pm\big(x_{\pm(i+1)}-x_{\pm i}\big)$,
$\eta_{\pm j}:=v_{\pm j}$, 
$i,j =1,2,\dots$,   are
completely independent, with exponential, respectively, normal
distributions: 
\bea
\notag
\prob{\xi_{\pm i}\in(x,x+\diff x)}=\ind_{\{x\ge0\}}e^{-x}\diff x,  \ \ 
\prob{\eta_{\pm j}\in(v,v+\diff
  v)}=\frac{1}{\sqrt{2\pi}}e^{-v^2/2}\diff v. 
\eea

We shall consider two different types of particle systems in this
paper.  Their state spaces will be
\bea
\notag
\OI
&:=&
\{(\op,\om,z,u,V): 
\opm\in\Opm, \ \ z\in\R_+, \ \ u\in[-1,1], \ \ V\in\R\},
\\[5pt]
\notag
\OII
&:=&
\{(\op,\om,z): 
\opm\in\Opm, \ \ z\in\R_+\}. 
\eea


We also define the natural projection between these spaces: 
\bea
\notag
\Pi:\OI\to\OII,
&\qquad&
\Pi(\op,\om,z,u,V):=(\op,\om,z). 
\eea

In order to define the relevant probability measures on the state
spaces $\OI$ and $\OII$ first we introduce some notation. 
Let the  random variables $W$ and $\zeta$ be independent and 
distributed as a standard Gaussian, respectively, as a standard
$\Gamma(2)$. Let $\gamma_2(z)$ be the density of the distribution of
$\zeta$, $\varrho(c)$ the density of the distribution of
$|W \zeta|$, and $\varphi_c(z)$ the density of the conditional
distribution of $\zeta$, given $|W \zeta|=c$:
\bea
\notag
\gamma_2(z)
&:=&
ze^{-z}, 
\\
\notag
\label{rhoc}
\varrho(c)
&:=&
\sqrt{\frac{2}{\pi}}\int_0^\infty \exp\{-z-\frac{c^2}{2z^2}\}\diff z, 
\\
\notag
\varphi_c(z)
&:=&
\frac{1}{\varrho(c)}\sqrt{\frac{2}{\pi}}\exp\{-z-\frac{c^2}{2z^2}\}.
\eea
Clearly, 
\bea
\notag
\gamma_2(z)= \int_0^\infty \varphi_c(z)\varrho(c) \diff c. 
\eea

The probability measures considered on the state  spaces $\OI$ and
$\OII$ are
$\mIM$ (defined on $\OI$) which 
depends on the positive parameter $M$, respctively, 
$\mIIc$ (defined on $\OII$) which
depends on the positive parameter $c$,  and  finally $\mII$  (also
defined on $\OII$) which is a mixture of the measures
  $\mIIc$: 
\bea
\label{Mmeasure}
\mIM(\diff \oI)
&:=&
\mp(\diff \op) \,\times\, \mm(\diff \om) \,\times\, 
\gamma_2(z) \diff z \,\times\, \frac12 \diff u \,\times\, 
\sqrt{\frac{M}{2\pi}} e^{MV^2/2}\diff V, 
\\
\label{cmeasure}
\mIIc(\diff \oII)
&:=&
\mp(\diff \op) \,\times\, \mm(\diff \om) \,\times\, 
\varphi_c(z)\diff z,
\\[5pt]
\notag
\mII(\diff \oII)
&:=&
\mp(\diff \op) \,\times\, \mm(\diff \om) \,\times\, \gamma_2(z) \diff z 
\\
\label{mixedmeasure}
&=& 
\int_0^\infty \mIIc(\diff \oII) \varrho (c) \diff c.
\eea
The measures  $\mIM$, respectively, $\mIIc$ will be the natural
Gibbs measures corresponding to the dynamics of our systems, to be
defined in the next subsection. 


\subsection{Dynamics}
\label{ss:dynamics}

We define the dynamics of the systems considered verbally, rather than
writing formulas.  
The 
two types of dynamics considered will be called
\emph{of type $\I$}, respectively,  \emph{of type $\II$}.
Their state spaces will be 
$\OI$, respectively,  $\OII$.
These  
will actually be families of dynamics parametrized by the fixed 
parameters $M>0$, respectively,  $c>0$.

\subsubsection{Dynamics of type $\I$:}
\label{sss:dyntypeI}

For precise formal definitions and basic facts about
these dynamics see \cite{sinaisoloveichik}, \cite{szasztoth1},
\cite{szasztoth2}.
The system consists of particles
indexed $\dots, -2,-1,0,+1,+2, \dots$. The system is observed from the
tagged particle of index  $0$. The tagged particle has mass $M$, the
other particles have unit mass. 
Positions and velocities of the particles in the system are encoded in
$(\op,\om,z,u,V)$ as follows: 
$V$ is the velocity of the tagged particle, 
$x_{\pm i}  \pm z(1\pm u)/2$ and 
$v_{\pm i}$ is the position, respectively, the velocity of the particle of
index $\pm i$, $i=1,2,\dots$.  
The untagged gas particles perform uniform motion on the line and 
don't interact between themselves, when two of them meet and cross
each other's trajectory they exchange their index. The tagged particle,
while isolated from the others, also performs uniform motion and
collides elastically at encounters with an untagged gas particle. 
At these collisions the outgoing velocities $V^{\mathrm{out}},
v^{\mathrm{out}}$  
are determined by the incoming velocities $V^{\mathrm{in}},
v^{\mathrm{in}}$   as follows: 
\bea
\label{coll}
V^{\mathrm{out}}=\frac{M-1}{M+1}V^{\mathrm{in}} + 
                 \frac{2}{M+1}v^{\mathrm{in}},
\qquad
v^{\mathrm{out}}=\frac{2M}{M+1}V^{\mathrm{in}} - 
                 \frac{M-1}{M+1}v^{\mathrm{in}}. 
\eea
Mind that the untagged gas particles never exchange their order with
the tagged particle and the index $\pm i$ of a particle denotes its
actual relative order with respect to the tagged particle.
 
The measure $\mIM$ defined in \eqref{Mmeasure} is Gibbs measure for
this dynamics, invariant 
for the system as seen from the tagged particle. It is a fact
(see \cite{sinaisoloveichik}, \cite{szasztoth1}, \cite{szasztoth2})
that the  
dynamics is $\mIM$-a.s. well defined: starting the system
distributed according to  $\mIM$, with probability 1 no multiple
collisions  will occur and the system remains locally finite
indefinitely. We denote by $\SIM_t$ the measure preserving
flow  defined by this dynamics on $(\OI, \mIM)$. 

The velocity and displacement  process of the tagged particle is 
\bea
\notag
&&
\VIM_t= \VIM_t(\oI):=V(\SIM_t \oI)
\\
\notag
&&
\QIM_t= \QIM_t(\oI):=\int_0^t \VIM_s(\oI)\diff s.
\eea

In section~\ref{s:Mtozero} it will be more convenient to describe 
the dynamics from a fixed exterior point of observation. The   
absolute locations of the gas particles in the system as seen from such a
fixed exterior frame of reference are 
\bea
\label{abspoz}
y_{\pm i}(t):=\QIM_t+x_{\pm i}(t),
\quad\text{where}\quad 
x_{\pm i}(t):=x_{\pm i}(\SIM_t \oI), 
\qquad
i=1,2,\dots.
\eea 

We also introduce the variables
\bea
\notag
&&
\tVIM_t= \tVIM_t(\oI):=\frac12\left(v_{-1}(\SIM_t \oI)+v_{+1}(\SIM_t
\oI)\right)
\\
\notag
&&
\tQIM_t= \tQIM_t(\oI):=\int_0^t \tVIM_s(\oI)\diff s. 
\eea
These are the velocity and position processes of the centre of mass of
the particles next to the right and to the left of the tagged
particle. We 
need the position process $\tQIM_t$ for later comparison with a
similar process defined for the dynamics of type $\II$ in the next
paragraph. Mind that the random process $t\mapsto \big(\QIM_t-\tQIM_t\big) $
is stationary and thus tight, uniformly for $t>0$. As a consequence,
$\big(\QIM_t/\sqrt{t}-\tQIM_t/\sqrt{t}\big)\to0$ in
$\mIM$-probability (actually 
 $\mIM$-a.s.) as $t\to\infty$. 

\subsubsection{Dynamics of type  $\II$:} 
\label{sss:dyntypeII}

The system consists of particles of unit mass
indexed $\dots, -2,-1,+1,+2, \dots$. Mind that there is no particle of
index $0$ in this system. The system is observed from the
centre of mass of particles of index $+1$ and $-1$, we call this the
central observation point. 
Positions and velocities of the
particles in the system are encoded in 
$(\op,\om,z)$ as follows: 
$x_{\pm i}  \pm z/2$,  
respectively,
$v_{\pm i}$ is the position relative to the central observation point,
respectively, the velocity of the particle of 
index $\pm i$, $i=1,2,\dots$. Clearly, $z$ denotes the distance
between the two central particles of index $+1$ and $-1$.   
Particles move uniformly on the line except for the two central
particles of index  
$+1$ and $-1$ which interact
via the inverse quadratic pair potential $U(z)$, or equivalently
repelling force $F(z)$:  
\bea
\label{pot}
U(z)=\frac{c^2}{2z^2}, 
\qquad
F(z)=\frac{c^2}{z^3},
\eea
where $c^2>0$ is a fixed parameter  and $z$ is the
distance between the two central particles. 
When two gas particles  meet and cross
each other's trajectory they exchange their index. But mind that due
to the strongly repulsive interaction between the two central
particles, these two will never meet and thus  particles will never
change the sign of their index. The 
index $\pm i$ of a particle denotes its 
actual relative order with respect to the central observation point. 

\medskip
\noindent
{\bf Remark.} 
In the literature of completely integrable Hamiltonian systems the
pair potential \eqref{pot} is usually called Calogero-Moser-Sutherland
interaction and leads to one of the most notorious completely
integrable 1-d systems, see \cite{calogero}, \cite{moser} and
\cite{sutherland} for the first original publications. 

\medskip

The measure $\mIIc$ defined in \eqref{cmeasure} 
is Gibbs measure for this dynamics, invariant
for the system as seen from the centre of mass of the two central
particles. It is again a fact that this 
dynamics is $\mIIc$-a.s. well defined.
We denote by $\SIIc_t$ the measure preserving
flow defined by this dynamics on $(\OII, \mIIc)$. 

The velocity and displacement process of the point of observation is
\bea
\notag
&&
\VIIc_t= \VIIc_t(\oII):=\frac12
\big(v_{-1}(\SIIc_t \oII) + 
     v_{+1}(\SIIc_t \oII) \big) 
\\
\notag
&&
\QIIc_t= \QIIc_t(\oII):=\int_0^t \VIIc_s(\oII)\diff s.
\eea
Again, absolute locations of the gas particles in the system as seen
from a fixed exterior frame of reference are expressed similarly to
\eqref{abspoz}.

\subsubsection{Stochastic processes considered}
\label{sss:processes}

In this paper we consider the following \emph{stochastic processes}:
\bea
\notag
\QIM_t
&=&
\QIM_t(\oI), 
\quad
\text{with random $\oI$ distributed according to}
\quad
\diff\mIM,
\\[3pt]
\notag
\tQIM_t
&=&
\tQIM_t(\oI), 
\quad
\text{with  random $\oI$ distributed according to}
\quad
\diff\mIM,
\\[3pt]
\notag
\QIIc_t
&=&
\QIIc_t(\oII), 
\quad
\text{with  random $\oII$ distributed according to}
\quad
\diff\mIIc,
\\[3pt]
\notag
\QII_t
&=&
\QIIc_t(\oII), 
\quad
\text{with random  $(\oII,c)$ distributed according to}
\quad
\diff\mIIc\varrho(c)\diff c,
\\[3pt]
\notag
&=&
\QIIc_t, 
\quad
\text{with random  $c$ distributed according to}
\quad
\varrho(c)\diff c.
\eea
This means that the process $\QII_t$ is a  $\varrho(c)\diff c$-mixture
of the processes  $\QIIc_t$

\section{Survey of earlier results}
\label{s:oldresults}

In this section we summarize the old results --
rigorously proved and numerical -- regarding
various limits for the motion of the tagged particle in the model of
type $I$. In Section~\ref{s:Mtozero} we formulate and prove a new result
concerning the $M\to 0$ behaviour of these systems.
In Section~\ref{s:newnumerics} we describe our new numerical results,
referring to this new model.

In all cases we are interested in the \emph{diffusive scaling limit of
the displacement of the tagged particle motion}, that is in the
asymptotics of the rescaled process 
\bea
\notag
t\mapsto A^{-1/2} Q_{At}, 
\quad\text{ as }\quad
A\to\infty.
\eea
Throughout the paper we denote by $A$ this scaling parameter. 

We shortly survey the existing results on the asymptotics of the
tagged particle motion in model of type $I$ in historical order. The
constants 
\bea
\notag
\usigma^2:=\sqrt{\pi/8}\approx0.627\dots,
\qquad
\osigma^2:=\sqrt{2/\pi}\approx0.798\dots 
\eea
will play a key role in the formulation of these results.

\subsection{The $M=1$ case:}
\label{ss:spitzer}

The case when the tagged particle has the same mass as the rest of the
gas particles was investigated and solved in Spitzer (1969)
\cite{spitzer}. For the roots of these ideas see also Harris (1965)
\cite{harris}. In \cite{spitzer} the following invariance principle is
proved: 
\bea
\notag
\text{for } M=1:
\qquad
A^{-1/2}\QIM_{At}\Rightarrow \osigma W_t,
\quad\text{ as }\quad 
A\to\infty, 
\eea
where $\Rightarrow$ stands for weak convergence of the sequence of
\emph{processes} (see \cite{billingsley} for weak convergence of
processes), and $W_t$ is a standard 1-d Borwnian motion. That is: 
$\osigma W_t$ is a Brownian motion of variance $\osigma^2$.

\subsection{The Ornstein-Uhlenbeck limit:}
\label{ss:holley}

Holley (1971) \cite{holley} considers the following limit when the
mass of the tagged particle is rescaled in the same order as the time
scale factor. Let $m\in(0,\infty)$ be fixed. Then
\bea
\notag
\text{for } M=mA:
\qquad
\big(A^{1/2}\VIM_{At}, A^{-1/2}\QIM_{At}\big)
\Rightarrow 
\big(\eta^m_t, \xi^m_t\big),
\quad\text{ as }\quad 
A\to\infty, 
\eea
where $\eta^m_t$ and $\xi^m_t$ are the Ornstein-Uhlenbeck velocity,
respectively, position processes defined by the SDEs
\bea
\notag 
\diff \eta^m_t=-\gamma(m)\eta^m_t \diff t +\sqrt{D(m)}\diff W_t,
\quad
\xi^m_t=\int_0^t\eta^m_s ds, 
\eea
with friction and dispertion parameters
\bea
\notag
\gamma(m):=\frac4m\sqrt{\frac2\pi},
\quad
D(m):=\frac8{m^2}\sqrt{\frac2\pi}.
\eea
For a version in higher dimensions of this type of result see D\"urr,
Goldstein, Lebowitz (1981) \cite{durrgoldsteinlebowitz}. 

It is important to remark (see \cite{szasztoth2}), that 
\bea
\notag
\xi^m_t\Rightarrow\usigma W_t, 
\quad\text{ as }\quad
m\to 0. 
\eea
This means that taking first Holley's limit, then $m\to0$ we obtain a
Wiener process of variance $\usigma^2$ as the diffusive scaling limit
of the displacement of the tagged particle.  

\subsection{Bounds for the limiting variance for any $M$:}
\label{ss:sinaisoloveichikszasztoth}

Sinai, Soloveichik (1986) \cite{sinaisoloveichik}, respectively,
Sz\'asz, T\'oth (1986) \cite{szasztoth1} consider the case of
arbitrary fixed  mass $M$ of the tagged particle. In these papers very
similar rersults are proved in completely different ways. The results
are summarized as follows: 
\bea
\notag
\text{for } M\ll A:
\qquad
\usigma^2 t \le 
\liminf_{A\to\infty} \var{A^{-1/2}\QIM_{At}} \le 
\limsup_{A\to\infty} \var{A^{-1/2}\QIM_{At}} \le 
\osigma^2 t.
\eea
Mind that these bounds are \emph{independent
of the mass of the tagged particle}. 
For surveys of these results see also \cite{ruellesinai},
\cite{soloveichik1}, \cite{soloveichik2}, \cite{szasztoth3}.

Any rigorous result regarding the
mass dependence of the limiting variance 
\bea
\notag 
\sigma_M^2:=\lim_{t\to\infty} \var{t^{-1/2}\QIM_{t}}
\eea
remains one of the most interesting open questions in this context
till today. The only known case is Spitzer's result
$\sigma_1=\osigma$. For numerical results see
subsection~\ref{ss:oldnumericalresults} below.

\subsection{Large mass Wiener limit:}
\label{ss:szasztoth2}

In order to interpolate between the $M=\text{const.}$ cases (see
subsection~\ref{ss:sinaisoloveichikszasztoth}) and 
Holley's limit (see subsection~\ref{ss:holley},) Sz\'asz and T\'oth
(1987) \cite{szasztoth2} considered the limit with
asymptotics $1\ll M \ll A$, as $A\to\infty$. Here the main result
is the following invariance principle: 
\bea
\label{szasztoth2}
\text{for } A^{1/2+\vareps}\ll M \ll A:
\qquad
A^{-1/2} \QIM_{At}\Rightarrow \usigma W_t,
\quad\text{ as }\quad
A\to\infty.
\eea
Actually, the scaling limit \eqref{szasztoth2} should hold for $1 \ll M
\ll A $ but the method of proof in  \cite{szasztoth2} based on a
coupling argument breaks down for  $1 \ll M \ll A^{1/2+\vareps}$.
For a survey of the results recalled in this and the previous
paragraph see also \cite{szasztoth3}.

\subsection{Earlier numerical results}
\label{ss:oldnumericalresults}

Following \cite{sinaisoloveichik} and \cite{szasztoth1} various
numerical investigations were performed in order to establish  the
mass dependence of the limiting variance: $M\mapsto\sigma^2_M$.  

\begin{figure}[hbt]
\psfrag{m}{$\sigma^{2,II,c}$}
\psfrag{sa}{$\sigmal^2$}
\psfrag{sf}{$\sigmau^2$}
\psfrag{sc}{$\sigma^2$}
\psfrag{c}{$c$}
\centering
\resizebox{12cm}{!}{\includegraphics[clip]{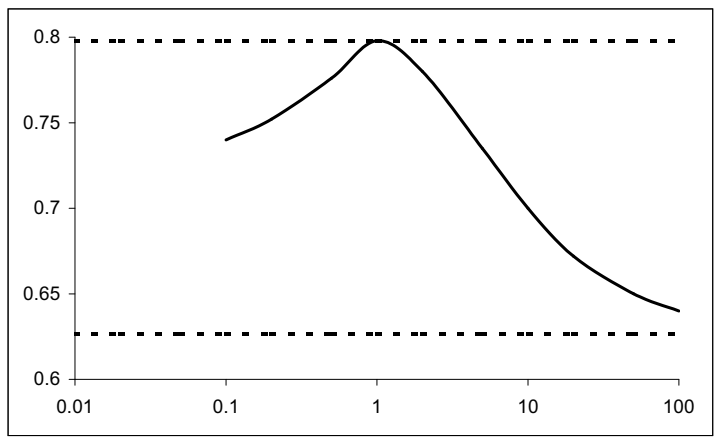}}
\caption{Qualitative dependence $M\mapsto\sigma^2_M$
suggested by earlier numerical works}
\label{fig:oldnumerics}
\end{figure}

The relevant numerical investigations performed  
in the late eighties, early
nineties are  published in Omerti, Ronchetti, D\"urr
(1986) \cite{omertironchettidurr}, Khazin (1987) \cite{khazin},
Boldrighini, Cosimi, Frigio (1990) \cite{boldrighinicosimifrigio},
Fernandez, Marro (1993) \cite{fernandezmarro}. 
These results clearly suggest the qualitative dependence
$M\mapsto\sigma^2_M$ shown in Figure~\ref{fig:oldnumerics}.

As about the $M\to0$ limit: in all these papers it is remarked that
the numerical simulations for small mass  of the tagged particle are
unreliable due to instability. On the other hand there was agreement
between all researchers  intrested in these questions that 
$\lim_{M\to0}\sigma_M^2 = \osigma^2\approx0.798\dots$
should hold. The ``starightforward argument'' was the following: the
tagged particle of extremely small mass must have 
very small effect on the system, vanishing as $M\to0$. So, in the
$M\to0$ limit the displacement of any marked  particle (in particular
the one next to the right of the tagged particle) will asymptotically
behave exactly like the tagged particle in Spitzer's  equal mass case,
cf. subsection~\ref{ss:spitzer}, above. So, the more recent and more
accurate 
numerical results published in Boldrighini, Frigio, Tognetti (2002)
\cite{boldrighinifrigiotognetti}, suggesting that 
\bea
\label{smallmasslimitnum}
\lim_{M\to0}\sigma_M^2 =: \sigma_0^2 \approx 0.74\dots
\eea
which is \emph{strictly} inbetween $\usigma^2\approx0.627\dots$ and
$\osigma^2\approx0.798\dots$, \emph{came as a surprise}. 

The results of the present note provide substantial theoretical
and independent numerical support of this surprising fact.

\section{The $M\to0$ limit of dynamics of type $\I$}
\label{s:Mtozero}

\begin{theorem}
\label{thm1}
Let $z\in\R^+$, $u\in[-1,1]$,  $W\in\R$ be fixed, $M_n\to0$,  
$V_n=M_n^{-1/2}W_n$ so that $W_n\to W\in\R$, and define
$c:=|Wz|$. Choose   
$\opm\in\Opm$ so that  
for all $n$ the dynamic trajectories 
${\cal S}^{\I,M_n}_t(\op,\om,z,u,V_n)$ and 
$\SIIc_t(\op,\om,z)$ are well defined for all
$t\in[0,\infty)$.  
(Mind that for any choice of $z$, $u$, $W$ and sequences $M_n$, $V_n$
  these  $\opm$-s are  of full $\mpm$ measure in
  $\Opm$.)  
Then for all $t\in[0,\infty)$ 
\bea
\notag
\lim_{n\to\infty}
\Pi{\cal S}^{\I,M_n}_t(\op,\om,z,u,V_n) = 
\SIIc_t(\op,\om,z). 
\eea
The convergence is uniform on compact intervals of time. 
\end{theorem}

\begin{proof}
Within this proof it is convenient to describe the systems of
particles as seen from a fixed external frame of reference: the
postion of the tagged particle (in the system of type $\I$) at time $t$
is $\QIM_t$, the postitions of the untagged gas particles are $y_{\pm
  i}(t)$, $i=1,2,\dots$, as given in \eqref{abspoz}. 

We have to prove that in the limit described in the theorem the
trajectories of the particles in system $I$ converge to the
coresponding trajectories in the limit system of type $\II$. 
Note that the particles with $i\ne \pm 1,0$ follow the same dynamical
rules in 
the two types of dynamics.  Thus we only need to understand how the
motion of the 
particles with indices $\pm 1$ can be approximated as $M\to 0$, if
they interact  
with the tagged particle according to the rules of
section~\ref{sss:dyntypeI}.   
  
As mentioned above, the particles with indices 
$-1$, $0$ and   $1$ have
positions $y_{-1}(t)\le Q_t \le  y_1(t)$  
and velocities $v_{-1}$, $V_t=W_t/\sqrt M$ and $v_1(t)$, respectively,
where $v_{\pm1}$ and $W$ are of order one.  
As $V$ is very large, the tagged particle  performs a full cycle: 
hits one of its neighbours, turns back, collides with the other
neighbour, and gets back to its  
initial position within a very short time $\diff t$. See
Fig~\ref{fig:flights} for insight. 

\begin{figure}[hbt]
\psfrag{dt}{$\diff t$}
\psfrag{t}{$t$}
\psfrag{y}{$y$}
\psfrag{m}{$-1.$}
\psfrag{0}{$0.$}
\psfrag{p}{$+1.$}
\centering
\resizebox{14cm}{!}{\includegraphics[clip]{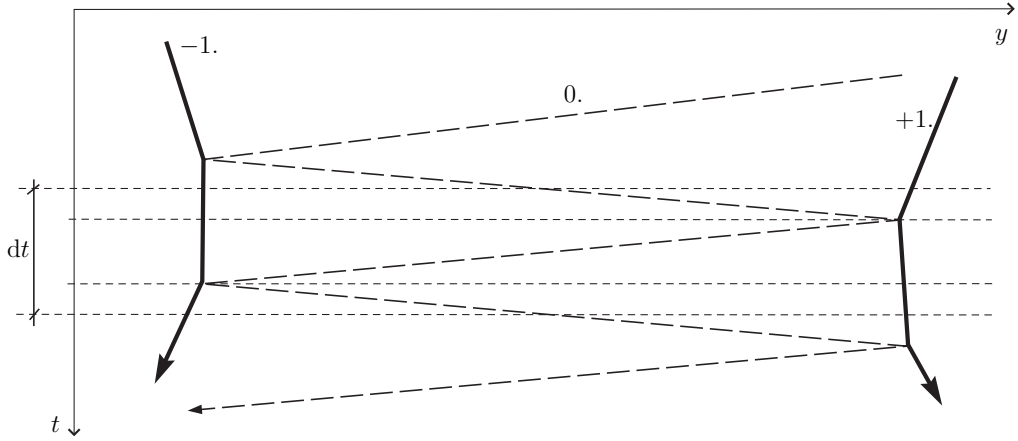}}
\caption{Successive collisions of the tagged particle of mass  $M\ll1$ with
its neighbours}
\label{fig:flights}
\end{figure}

To investigate how the system evolves in this time interval $\diff t$,
two  successive collisions should be taken into account. 
We may assume that $W>0$ (the case of negative $W$ is analogous), 
and thus the tagged particle collides first with the particle of index 
$1$,  and then with that of index $-1$. These collisions split the
time interval  into three smaller subintervals. 
Let us expand the formulas of \eqref{coll} in the limit as $M\to
0$. This way we may calculate    the velocities of the three particles
in the three subintervals, see  Table~\ref{tab:collapp}.  
Mind that the order of magnitude of the incoming velocities are as
follows
\bea
\notag
V=W/\sqrt M\, \asymp \,M^{-1/2},
\qquad
v_{\pm1}\, \asymp \,1. 
\eea
\begin{table}[h!bt]
\begin{tabular}{l||c|c|c}
\hline
&&&
\\[-5pt]
           &
1. interval&
2. interval&
3. interval
\\[-5pt]
&&&
\\
  \hline 
  \hline
&&&
\\[-5pt]
  particle 1&
  $v_1$&
  $v_1+2M V+{\cal O}(M)$&
  $v_1+2M V+{\cal O}(M)$
\\[-5pt]
&&&
\\[-5pt]
  \hline
&&&
\\[-5pt]
  particle 0&
  $V$&
  $-V+2v_1+{\cal O}(\sqrt{M})$&
  $V-2v_1+2v_{-1}+{\cal O}(\sqrt{M})$
\\[-5pt]  
&&&
\\[-5pt] 
  \hline
&&&
\\[-5pt] 
  particle -1& 
  $v_{-1}$&
  $v_{-1}$&
  $v_{-1}-2MV+{\cal O}(\sqrt{M})$
\\[-5pt]   
&&&
\\   
\hline
\end{tabular}
\caption{Velocities in the $M\to 0$ approximation}
\label{tab:collapp} 
\end{table}

In particular, the tagged particle 
reverts its velocity at collisions, thus its speed is constant, more
precisely, equal to  
$|V|+\cO(1)=|W/\sqrt{M}|+\cO(1)$ 
throughout the investigated time interval.
This implies $dt\asymp \sqrt{M}$. 
Furthermore, the velocities $v_1$ and $v_{-1}$  
remain $\cO(1)$, and thus, the particles of index $1$ and $-1$ remain
$\cO(\sqrt{M})$-close 
to their original positions  $y_1$ and  $y_{-1}$ in the investigated 
interval.   
In accordance with the notations of subsection~\ref{ss:statespace}, let
us introduce 
$z=y_1-y_{-1}$ for brevity. By the above observations the distance of
the two non-tagged particles  
remains $\cO(\sqrt{M})$-close to $z$. Now we may calculate the leading
term in $\diff t$ : 
\bea
\label{dt}
\diff t=\frac{2z + \cO(\sqrt{M}) }{V +\cO(1)} = \frac{2z}{W}\sqrt{M} + 
\cO(M). 
\eea
Let us denote the amount with which the velocities change at the
collisions  
(and thus, during the studied time interval) by $\diff v_1$ and $\diff
v_{-1}$.   

     
Referring to Table~\ref{tab:collapp} we get
\bea
\label{dv}
\diff v_1=2W\sqrt{M}+\cO ({M}), 
\qquad 
\diff v_{-1}=-2W\sqrt{M}+ \cO({M}). 
\eea
 
Formulas \eqref{dt} and \eqref{dv} altogether imply that, 
as $M\to 0$, $v_1(t)$ and $v_{-1}(t)$ approach (piecewise)
differentiable functions, and 
\bea
\label{vdot}
\dot{v}_1=\frac{W^2}{z}, \qquad \dot{v}_{-1}=-\frac{W^2}{z}.
\eea

Referring again to Table~\ref{tab:collapp}, we may calculate the
amount of 
change in the velocity of the tagged particle during the time interval
$\diff t$. 
We get
\bea
\label{dv0}
\diff V=2v_{-1}-2v_1 +\cO(\sqrt{M}),
\qquad {\rm thus} \qquad 
\diff W=\sqrt{M}(2v_{-1}-2v_1) +\cO({M}).
\eea
Formulas \eqref{dt} and \eqref{dv0} altogether imply
\bea 
\notag
\label{dvdW}
\dot{W}=-\frac{(v_1-v_{-1})W}{z}=-\frac{\dot{z} W}{z},
\eea
where we have used that $v_1-v_{-1}=\dot{y}_1-\dot{y}_{-1}=\dot{z}$.

Integrating this differential equation we find that $|W z|$ is an
integral of motion. This means that in the $M\to 0$ limit,
$c=|W(t) z(t)|$ is \emph{constant}
during a time interval when the tagged particle has the same neighbours,
and within such an inverval, \eqref{vdot} gives that
the positions and velocities for the particles with indices $\pm 1$
satisfy the coupled differential 
equations
\bea
\notag
\dot{x}_1=v_1, \qquad & \dot{v}_1={c^2}/{z^3}, 
\\[5pt]
\notag
\dot{x}_{-1}=v_{-1}, \qquad & \dot{v}_{-1}=-{c^2}/{z^3}. 
\eea
This is in agreement with the formula \eqref{pot} for the potential
describing the systems of type $\II$.

Notice now that the value of $c=|W(t) z(t)|$
also remains constant when one of the neighbours of the tagged particle
`meets' another gas particle, and the neighbour is replaced by that
new particle. 
At such a time moment the values of both $W$ and $z$ are unchanged.
Thus for any $t>0$ $|W(t) z(t)|=|W(0) z(0)|=c$
according to the choice of $c$ in the formulation of Theorem~\ref{thm1}.
This completes the proof of Theorem~\ref{thm1}. 
\end{proof}


\noindent
{\bf Remark.} 
Note that taking the $c\to 0$ limit of the dynamics of
type $\II$, 
we recover the dynamics of type $\I$ with equal masses: 
the interaction between the two central particles becomes hard core
specular collision. So, in  this double limit ($M\to0$ and than
$c\to0$)) the 
system behaves indeed as Spitzer's model, see
subsection~\ref{ss:spitzer}.  

\medskip

Recall that according to \eqref{mixedmeasure} the measure $\mII$
which is the projection of the measures $\mIM$ on the state space $\OII$,
is the $\varrho(c)\diff c$-mixture of the Gibbs measures $\mIIc$. This
implies that Theorem~\ref{thm1} has the following immediate corollary:

\begin{corollary*}
\label{cor}
Let $M\to 0$. 
For any fixed $0<T<\infty$, the sequence of processes
$[0,T]\ni t\mapsto\tQIM_t$ converges weakly (in distribution) to the
process to $[0,T]\ni t\mapsto\QII_t$.  
\end{corollary*}


\section{Numerical results on systems of type $\II$}
\label{s:newnumerics}

\subsection{Generalities}
\label{ss:generalities}

In this section we describe numerical investigations aimed at calculating the limiting variance
\bea
\notag
\sigma^2:=\lim_{t\to\infty}t^{-1}\var{Q_t}
\eea
for the systems of type $\II$. We shall also comment on how these results are
related to the $M\to 0$ limit of the variance for the systems of type $\I$, as established 
numerically in \cite{boldrighinifrigiotognetti}.

These simulations of the systems of type $\II$ 
were done by following a number of particles for some fixed time $T$. 
The particles followed were those who were less than $10 T$ far away
from the point of 
observation in the beginning. It is easy to check that with this method the
probability of not following a particle that would indeed participate in the
interaction, is negligible in all the cases we looked at.

Numerial simulation of the dymanics of type $\II$ is relatively fast, since the
equation of motion for the two particles interacting via the potential
\eqref{pot} can be solved explicitly. (This observation is at the heart
of the complete solvability of the Calogero-Moser-Sutherland model.)

The simulation for time $T$ was repeated over a sample of $N$ initial
conditions 
chosen independently according to the appropriate stationary Gibbs
distribution. 
From this sample, the empirical variance was calculated for
$\var{Q_t}$ as a function of $t$.

\begin{figure}[hbt]
\psfrag{m}{$\sigma^{2,II,c}$}
\psfrag{sa}{$\sigmal^2$}
\psfrag{sf}{$\sigmau^2$}
\psfrag{sc}{$\sigma^2$}
\psfrag{c}{$c$}
\centering
\resizebox{12cm}{!}{
\includegraphics[clip]{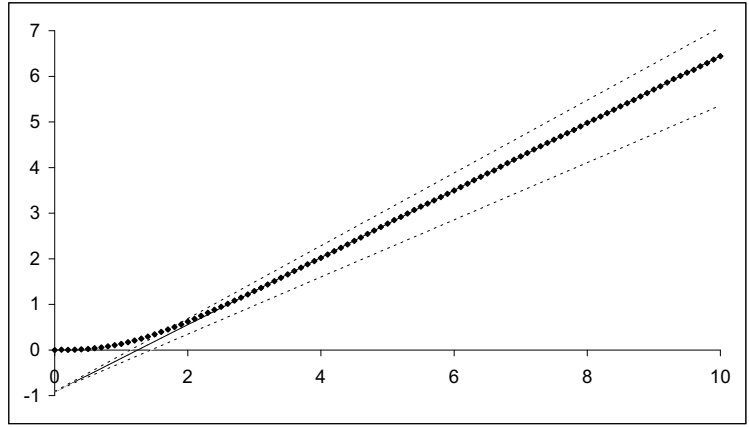} }
\caption{$\var{Q_t}$ as a function of $t$ in a typical simulation}
\label{fig:alap_illesztes}
\end{figure}

The result of a typical simulation can be seen in
Figure~\ref{fig:alap_illesztes}. 
The solid line is the best linear fit for the tail, while the dashed lines have
slope $\sigmal^2$ and $\sigmau^2$, and are drawn for comparison.

As we can see, $\var{Q_t}$ does appear to be asymptotically linear.
To read the limit $\displaystyle \lim_{t\to\infty}t^{-1}\var{Q_t}$
from the graph, we needn't perform such a long simulation, during
which this limit is well approached: the slope of the asymptote
can be found with a good accuracy much sooner. Thus all the limits
given in the paper are obtained using this technique, and the time
interval for the simulation is typically between $T=10$ and
$T=50$. In exchange, the size of the sample can be very big --
actually, samples up to $N=10^{7}$ were used.

Finally, the statistical error of the calculated values was estimated
by simply 
repeating the whole procedure about $20$ times and calculating the
standard deviation of the values obtained.

A detailed description of the numerical simulation and the source code for
the applied program can be found in \cite{program}.

\subsection{The systems with fixed $c$}
\label{ss:fixedcnumerics}

\begin{figure}[hbt]
\psfrag{m}{$\sigma^{2,II,c}$}
\psfrag{sa}{$\sigmal^2$}
\psfrag{sf}{$\sigmau^2$}
\psfrag{sc}{$\sigma^2$}
\psfrag{c}{$c$}
\centering
\resizebox{12cm}{!}{
\includegraphics[clip]{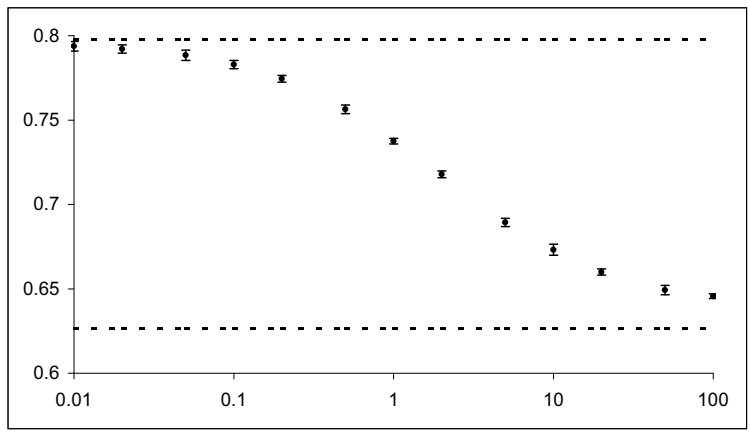} }
\caption{$c$-dependence of the limiting variance for systems of type $II$}
\label{fig:c-fugges}
\end{figure}

We simulated numerically the  dynamics of type $\II$ for various
fixed values of the parameter $c$ ranging between 0.01 and 100. We
started the system from samples of the stationary Gibbs distribution 
$\mIIc$ and computed the limiting variance
\bea
\notag
\sigma^2_c:=
\lim_{t\to\infty}
t^{-1}\var{\QIIc_t}. 
\eea
We found the  $c\mapsto \sigma^2_c$ dependence of the limiting
variance  as
shown in Figure~\ref{fig:c-fugges}. 
We see that $\lim_{c\to0}\sigma^2_c \to \sigmau^2$, 
which is no surprise, since in the $c\to 0$ limit 
the system indeed
behaves like the system of type $\I$ with $M=1$, which is known to
have $\sigma^2_{M=1}=\sigmau^2$. See subsection~\ref{ss:spitzer} and the
Remark after the proof of Theorem 1.




On the other hand, it is interesting to see that as $c\to\infty$,
the limiting variance decreases, and even seems to approach
a value near the lower limit $\sigmal^2$, but not quite reaching this
lower bound. We plan to return to this phenomenon in the forthcoming
paper \cite{balinttothtoth2}.

\subsection{The mixed system}
\label{ss:mixednumerics}

We computed the numerical value of the limiting variance for the
mixture of dynamics  in two different ways. 

First, we computed numerically the value of 
\bea
\notag
\sigma^2_{\text{mix,1}}:=\lim_{t\to\infty} 
t^{-1}\var{\QII_t}. 
\eea
We have done it in the following way: 
we sampled the initial conditions $\oII$ according to the distribution
$\mII$ and, independently, a standard normal variable $W$. Then we
computed $c:=|W  z|$, where $z$ was the distance between the two
central particles in the initial configuration $\oII$. This
\emph{random} value $c$ served as strength parameter in the
interaction potential \eqref{pot}, with which the dynamics $\SIIc_t$ was 
computed. 

Second, using the data obtained for $\sigma_c^2$ in the fixed $c$
computations (see subsection~\ref{ss:fixedcnumerics}), we computed the
mixture  
\bea
\notag
\sigma^2_{\text{mix,2}}:=\int_0^{\infty} \sigma^2_c \rho(c) \diff c,
\eea 
which, of course, in principle must give the same value as the
previous computation. 

Indeed, in the two cases we obtained the numerical values
$$
\sigma^2_{\text{mix,1}}= 0.736 \pm 0.003, \quad
\sigma^2_{\text{mix,2}}= 0.740 \pm 0.003.
$$
This result is very interesting, since it coincides exactly (well within
statistical error) 
with the $M\to 0$ limit of the variance calculated numerically in
\cite{boldrighinifrigiotognetti}, see \eqref{smallmasslimitnum}. This
means that there is indeed  \emph{continuity} in the limiting variance
as $M\to 0$. 

We remark that the function $c\mapsto\sigma^2_c$ shown in
Figure~\ref{fig:c-fugges} can be fit with amazing accuracy by the function of
the simple form
$\sigma^2_c=\frac{A_1-A_2}{1+\left(\frac{c}{c_0}\right)^p}+A_2$, where
$A_1=0.796$, $A_2=0.638$, $c_0=1.981$, $p=0.792$.

\vskip2cm
\noindent
{\bf Acknowledgement.} 
We thank Henk van Beijeren for pointing out the relation of our result
with the Calogero-Moser-Sutherland interaction. 
The research work of
the authors is partially supported by the following OTKA (Hungarian
National Research Fund) grants:  
F 60206 (for P.B.),
K 60708 (for B.T.), 
T 46187 (for P.T.)
and 
TS 49835 (for P.B. and P.T.); and by the Bolyai scholarship of the
MTA (for P.B.).

\vskip2cm

\hbox{
\phantom{M}
\hskip7cm
\vbox{\hsize8cm
{\sc
\noindent
Address of authors: 
\\[8pt]
\dots-P\'eter-B\'alint-T\'oth-P\'eter-\dots
\\
Institute of Mathematics
\\
Technical University of Budapest
\\
Egry J\'ozsef u. 1
\\
H-1111 Budapest, Hungary}
\\[5pt]
e-mail: 
{\tt \{pet,balint,mogy\}{@}math.bme.hu}
}}

\end{document}